\documentclass[12pt]{amsart}
\usepackage{amsmath,amsfonts,amssymb}

\def\ssum{\mathop{\sum\!\sum}}

\def\le{\leqslant}
\def\ge{\geqslant}

\renewcommand{\mod}{\mathop{\rm{mod}}}

\numberwithin{equation}{section}

\newtheorem{cor}{COROLLARY}[section]
\newtheorem{prop}[cor]{PROPOSITION}

\theoremstyle{definition}

\newtheorem{theorem}{THEOREM}

\begin{document}
\title{\bf Twin Primes Via Exceptional Characters}
\author[Friedlander]{J.B. Friedlander$^*$}
\thanks{$^*$\ Supported in part by NSERC grant A5123}
\author[Iwaniec]{H. Iwaniec$^{**}$}
\thanks{$^{**}$\ Supported in part by NSF grant DMS-1406981}

\maketitle

{\bf Abstract:} 
We give an account of the arguments that lead from the assumption 
of the existence of exceptional characters to the asymptotics in 
related ranges for the counting function of twin primes. 

\section{\bf Introduction}

The twin prime conjecture asserts that there are infinitely many 
primes $p$ such that $p+2$ is also prime. More generally, 
given a positive even number $h$, we expect the asymptotic formula
\begin{equation}\label{eq:1.1}
S_h(x) = \sum_{n\le x} \Lambda(n)\Lambda(n+h)\sim BC(h)x 
\end{equation} 
to hold as $x\rightarrow \infty$, where $B$ is the absolute constant 
\begin{equation}\label{eq:1.2}
B=2\prod_{p>2}\Bigl(1-\frac{1}{p-1}\Bigr)\Bigl(1-\frac{1}{p}\Bigr)^{-1}
\end{equation} 
and $C(h)$ depends mildly on $h$, namely
\begin{equation}\label{eq:1.3}
C(h) = \prod\limits_{\substack{p|h\\p>2}} \, 
\Bigl(1 - \frac{1}{p-1}\Bigr)^{-1}\ . 
\end{equation} 
This is a conjecture of Hardy and Littlewood. 
Actually, it is easy to predict more general formulas, 
such as (3.7) of [FI1] by exploiting the assumption of 
randomness of the M\"obius function in conjunction with sieve ideas. 

In these notes we present elementary arguments of sieve type which 
yield the asymptotic formula (1.1)  with an error term estimated by means 
of the series 
\begin{equation}\label{eq:1.4}
L(1,\chi) =\sum_1^{\infty}\chi(n)n^{-1} 
\end{equation}
with a primitive, real character $\chi (\mod D)$. 

\begin{theorem}  
Let $x\ge D^{3500}$. For any even positive number $h$ we have 
\begin{equation}\label{eq:1.5}
S_h(x)=BC(h)x +O\bigl(L(1,\chi)x \log x +x/\log x\bigr) 
\end{equation}
where the implied constant depends only on $h$. 
\end{theorem} 

The result is unconditional but is meaningful only if $L(1,\chi)$ 
is sufficiently small. We put 
\begin{equation}\label{eq:1.6}
\eta(D) = L(1,\chi) \log D\ . 
\end{equation}
If $\eta (D)$ is small we call the character $\chi$ exceptional. 
For such characters the error term in (1.5), say $xE(x)$, is  
also relatively small, namely
\begin{equation}\label{eq:1.7}
E(x) \ll \eta(D)A +1/\log D, 
\quad {\rm if} \quad D^{3500} \le x \le D^A\ . 
\end{equation}
Hence, we conclude the following result of Heath-Brown [H-B]. 
\begin{cor}  
If there are infinitely many exceptional characters then there are 
infinitely many twin prime numbers. 
\end{cor}

{\bf Acknowledgement:} The second author thanks 
ETH-ITS in Zurich for wonderful working 
conditions and financial support for his visit of June-July 2016, 
during which time this paper was completed. 

\section {\bf A Partition of $S_h(x)$}

First, for notational simplicity, we write 

$$
S_h(x) = \ssum_{m-n=h}\Lambda(m)\Lambda(n) + O\bigl(h(\log x)^2\bigr) 
$$
where here and in the following we understand, but do not display, the 
conditions $0<m, n \le x$, $(mn,h)=1$. Next, into 
$\Lambda(m)$, $ \Lambda(n)$ we introduce upper-bound sieve factors 
$\theta(m)$, $ \theta(n)$ which are almost redundant. Specifically, 
let $(\xi_q)$ be sieve weights of level $y$ and range $P(z)$ which is the 
product of all primes $p<z$, $p$ not dividing $h$. 
This means we have real numbers 
$\xi_q$ with $\xi_1 =1$, $|\xi_q| \le 1$, and $\xi_q= 0$ unless $q\mid P(z)$, 
$q<y$, and such that for any positive integer $m$ 
\begin{equation}\label{eq:2.1}
\theta (m) =\sum_{q\mid m}\xi_q \ge 0\ . 
\end{equation}
Note that $\theta(m)$ is bounded by the divisor function $\tau(m)$. We get 
\begin{equation}\label{eq:2.2}
S_h(x) = \ssum_{m-n=h}\theta(m)\theta(n)\Lambda(m)\Lambda(n) 
+ O\bigl((h+z)(\log x)^2\bigr)\ . 
\end{equation}
The level $y$ and range $z$ of the sieve will be specified later. 

Now we are ready to decompose $\Lambda$ in terms of the following 
Dirichlet convolutions: 
\begin{equation}\label{eq:2.3}
\lambda = \chi * 1\ , \quad \lambda' =\chi * \log = \lambda * \Lambda\ , 
\quad \nu=\mu \chi *\mu \ . 
\end{equation} 
Note that $0\le \Lambda \le \lambda' \le \tau \log$ 
and $|\nu| \le \lambda \le \tau$. By M\"obius inversion 
\begin{equation}\label{eq:2.4}
\Lambda = \lambda' *\nu \ . 
\end{equation} 
Having in mind that $\nu$ can be lacunary, we split $\Lambda = \Lambda^* 
+ \Lambda_*$ with 
\begin{equation}\label{eq:2.5}
\Lambda^*(n) = \sum_{\substack{ab=n\\b<y}}\lambda'(a)\nu(b)\ , \quad\quad 
\Lambda_*(n) = \sum_{\substack{ab=n\\b\ge y}}\lambda'(a)\nu(b) \ . 
\end{equation} 
Note that, for economy of notation, we have chosen the splitting 
parameter $y$ to be the same as that for the sieve level. The 
splitting parameter will be relatively small ($y\le x^{1/9}$), 
so the variable $b$ in $\Lambda^*(n)$ is quite short and $\Lambda^*(n)$ 
looks like the divisor-type function $\lambda' * \log $ which can be 
analyzed by various means (such as by Dirichlet's switching-divisors 
technique, the $\delta$-method, the circle method, 
the spectral decomposition of Poincar\'e series). On the other hand, 
the variable $b$ in $\Lambda_*(n)$ is sufficiently long that the 
lacunarity of the factor $\nu(b)$ will kick in to yield crude but 
admissible estimates. 

Writing 
$ \Lambda\Lambda = \Lambda^*\Lambda^* +\frac12 (\Lambda + \Lambda^*)\Lambda_*   
+\frac12 \Lambda_* (\Lambda + \Lambda^*)$, 
we get 
\begin{equation}\label{eq:2.6}
S_h(x) = S_h^*(x) +\tfrac12 T_h(x) +\tfrac12 T_{-h}(x)
+ O\bigl((h+z)(\log x)^2\bigr)\ , 
\end{equation}
where 
\begin{equation}\label{eq:2.7}
S_h^*(x) = \ssum_{m-n=h}\theta(m)\theta(n)\Lambda^*(m)\Lambda^*(n)\ ,
\end{equation}
\begin{equation}\label{eq:2.8}
T_h(x) = \ssum_{m-n=h}\theta(m)\theta(n)\bigl(\Lambda(m) +\Lambda^*(m)\bigr) 
\Lambda_*(n)\ . 
\end{equation}

\section {\bf Estimates for Divisor-Like Functions}

In $T_h(x)$ we estimate 
$|\Lambda_*(n)| \le a(n)\log x$, $|\Lambda^*(m)| \le b(m)\log x$ with
\begin{equation}\label{eq:3.1}
a(n)= \sum_{\substack{ab=n\\b\ge y}}\tau(a)\lambda(b)\ , \quad 
b(m)= \sum_{\substack{ab=m\\b< y}}\tau(a)\tau(b)\ .
\end{equation}
Here we have (note that $2\alpha^{13}\le \tau_{14}(p^{\alpha})$ in Proposition 
22.12 of [FI1]), 
$$
\tau(a)\le 2^{31}\sum_{\substack{c|a\\c\le a^{1/4}}}\tau_{14}(c)\ .
$$
Hence
$$
b(m)\le 2^{31}\sum_{\substack{d|m\\d<m^{1/4}y^{3/4}}}\tau_{16}(d)\ .
$$
For $m\le x$ and $y\le x^{1/9}$ we have divisors $d<x^{1/3}$. Moreover, we have 
$$
\tau_{16}(d)\le 4\tau_{17}(d)\varphi(d)/d 
$$
because $\tau_{17}(p^{\alpha})/\tau_{16}(p^{\alpha})
= 1 +\alpha/16 \ge 17/16\ge (1-1/p)^{-1}$ if $p>13$. 

Finally, we estimate $\tau_{17}(d)\varphi(d)/d$ by 
\begin{equation}\label{eq:3.2}
\gamma(d) =\frac{\varphi(d)}{d} \sum_{d_1\ldots d_r=d} d_1^{\varepsilon_1}\ldots 
d_r^{\varepsilon_r}\ , 
\end{equation}
where $\varepsilon_1\ldots \varepsilon_r$ are distinct small positive numbers 
with $r =17$. We obtain $b(m)\le2^{33}c(m)$ with 
\begin{equation}\label{eq:3.3}
c(m)= \sum_{\substack{d|m\\d<x^{1/3}}}\gamma(d)\ . 
\end{equation}
From the above estimates we conclude that
\begin{equation}\label{eq:3.4}
|T_h(x)| \le 2^{34}V_h(x)(\log x)^2\ , 
\end{equation}
with
\begin{equation}\label{eq:3.5}
V_h(x)=  \ssum_{m-n=h}\theta(m)\theta(n)c(m)a(n)\ . 
\end{equation}

The reason for making a slight deformation of the divisor function 
$\tau_{16}(d)$ in the above argument is so as to reach a simple generating 
Dirichlet series, specifically 
\begin{equation}\label{eq:3.6}
D(s)=\sum_d\frac{\gamma(d)}{\varphi(d)} d^{-s} =\zeta(s+1-\varepsilon_1)
\ldots\zeta(s+1-\varepsilon_r) \ .  
\end{equation}
We shall take $\varepsilon_i =i/\log z$ for $1\le i \le r = 17$. 

\section {\bf Estimation of $V_h(x)$} 

By (2.1), (3.3) and (3.5) we arrange $V_h(x)$ into sums of congruence sums 
\begin{equation}\label{eq:4.1}
V_h(x)= \sum_{\substack{q<y\\(q,h)=1}}\xi_q\sum_{\substack{d<x^{1/3}\\(d,h)=1}} \gamma(d) 
\sum_{\substack{n\equiv -h ([d,q])\\(n,h)=1}}\theta(n)a(n)
\ , 
\end{equation}
where $n$ runs over the segment $\max (0,-h)<n\le \min (x,x-h)$. We 
extend this segment to $0<n\le x$ up to an error term of size 
$O(|h|(\log x)^{20})$. Note that the sequence ${\mathcal A} = (\theta(n)a(n))$ 
with $(n,h)=1$, $0<n\le x$ admits a level of distribution 
$x^{\frac 12 -\varepsilon}$. This can be achieved by the large sieve method. 
Therefore, our congruence sums are equidistributed over reduced residue classes 
apart from an error term $O(x^{1-\delta})$ provided that $Dy^2\le x^{1/9}$. We 
have 
\begin{equation}\label{eq:4.2}
V_h(x)= M_h(x) +O\bigl(|h| (\log x)^{20} +x^{1-\delta} \bigr) 
\, 
\end{equation}
with
\begin{equation}\label{eq:4.3}
M_h(x)= \sum_{\substack{q<y\\(q,h)=1}}\xi_q\sum_{\substack{d<x^{1/3}\\(d,h)=1}} 
\frac{\gamma(d)}{\varphi([d,q])} 
\sum_{\substack{0<n\le x\\(n,hdq)=1}}\theta(n)a(n)\ , 
\end{equation}
where $\delta>0$ and the implied constant is absolute. 

We are going to execute the summation over $d$ first, which for given 
$n$ and $q$, is equal to 
\begin{equation}\label{eq:4.4}
L(x) = \sum_{\substack{d<x^{1/3}\\(d,hn)=1}}\frac{\gamma_d}{\varphi\bigl([d,q]\bigr)}
 = \frac{1}{\varphi(q)}\sum_{c|q^{\infty}}\gamma_c\frac{(c,q)}{c} 
\sum_{\substack{d<x^{1/3}/c\\(d,hnq)=1}} 
\frac{\gamma(d)}{\varphi(d)}\ . 
\end{equation}
Note that the divisors $c$ of $q^{\infty}$ with $c>x^{1/9}$ conribute 
a negligible amount. For smaller $c$ we evaluate the sum over $d$ by 
contour integration of the series 
$$
D_v(s) =\sum_{(d,v)=1} \frac{\gamma(d)}{\varphi(d)} = P_v(s)D(s)\ ,
$$
where $D(s)=D_1(s)$ is given by the product of zeta-functions (3.6) and 
$P_v(s)$ removes the local factors of $D(s)$ at primes $p|v$; 
\begin{equation}\label{eq:4.5}
P_v(s)=\prod_{p|v}\bigl(1-p^{\varepsilon_1-s-1}\bigr)\ldots 
\bigl(1-p^{\varepsilon_r-s-1}\bigr)\ . 
\end{equation}
Hence $D_v(s)/s$ has simple poles at $s=\varepsilon_0 =0$ 
and at $s=\varepsilon_1\ ,\ldots \varepsilon_r$ with residue
$$
R_v(i) = R(i)P_v(\varepsilon_i)
$$
where $R(i)$ is the residue of $D(s)/s$ at $s=\varepsilon_i$, that is 
$$
R(i) =\kappa(i) \prod^r_{\substack{j=0\\j\neq i}}
\zeta(1+\varepsilon_i -\varepsilon_j)
$$
with $\kappa(0)=1$ and $\kappa(i)=1/\varepsilon_i$ if $1\le i \le r$. 
Note that 
\begin{equation}\label{eq:4.6}
R(i) \asymp (\log z)^r\ . 
\end{equation}
By complex integration, the inner sum over $(d,hnq)=1$ in (4.4) is 
equal to 
\begin{equation}\label{eq:4.7}
\sum^r_{i=0}R(i) P_{hnq}(\varepsilon_i)\bigl(x^{1/3}/c\bigr)^{\varepsilon_i}
\end{equation} 
up to a small, negligible error term. Next we compute the resulting 
multiplicative functions (the sieve density) 
\begin{equation}\label{eq:4.8}
g_ {\varepsilon}(q) = \frac{P_q(\varepsilon)}{\varphi(q)}
\sum_{c|q^{\infty}}\gamma(c)(c,q)c^{-\varepsilon-1}
\end{equation} 
for every $\varepsilon =\varepsilon_i$, $0\le i \le r$. At primes we have 
$$
g_ {\varepsilon}(p) = \frac{P_p(\varepsilon)}{p-1}
\Bigl(1+p\sum^{\infty}_{\alpha =1}\gamma(p^{\alpha})p^{-\alpha(\varepsilon+1)}\Bigr)
=1-\frac{p-2}{p-1}P_p(\varepsilon)\ . 
$$
Note that $0<P_p(\varepsilon)\le 1$, so $0<g_{\varepsilon}(p)<1$. More 
precisely, we have 
$$
P_p(\varepsilon) = 1 - \sum^r_{j=1}p^{\varepsilon_j - \varepsilon -1} 
+O\bigl(p^{-2}\bigr)\ , 
$$
so 
\begin{equation}\label{eq:4.9}
g_{\varepsilon}(p) = p^{-1} + \sum^r_{j=1}p^{\varepsilon_j - \varepsilon -1} 
+O\bigl(p^{-2}\bigr)\ , 
\end{equation} 
for every $\varepsilon = \varepsilon_i$, $0\le i\le r$. 
Since the $\varepsilon_i =i/\log z$ are small we have essentially 
a sieve problem of dimension $r+1$. 

By the above computations, we get 
$$
L(x) = \sum^r_{i=0}R(i)P_{hn}(\varepsilon_i)
g_{\varepsilon_i}(q)x^{\varepsilon_i /3}\ ,  
$$
up to a small, negligible error term. Hence, (4.3) becomes 
\begin{equation}\label{eq:4.10}
M_h(x)= \sum^r_{i=0}R(i)P_{h}(\varepsilon_i)x^{\varepsilon_i /3}
\sum_{\substack{n\le x\\(n,h)=1}}\theta(n)a(n)P_n(\varepsilon_i)G_{hn}(\varepsilon_i) 
\end{equation} 
where
\begin{equation}\label{eq:4.11}
G_v(\varepsilon) = \sum_{\substack{q<y\\(q,v)=1}}\xi_qg_{\varepsilon}(q) 
\end{equation} 
up to a small, negligible error term. By sieve methods, 
$$
G_v(\varepsilon) \asymp \prod_{\substack{p<z\\ p \nmid v}}
\bigl(1-g_{\varepsilon}(p)\bigr) 
\asymp \prod_{\substack{p<z\\ p \nmid v}}\bigl(1-p^{-1}\bigr)
\bigl(1-p^{\varepsilon -\varepsilon_1-1}\bigr)\ldots 
\bigl(1-p^{\varepsilon -\varepsilon_r-1}\bigr)\ ,
$$
provided $z^{4r} \le y$. Actually, this condition is not necessary if 
we claimed only the upper bound for $G_v(\varepsilon)$. 
For $s=1+O(1/\log z)$ we have 
$$
\prod_{p<z}\bigl(1-p^{-s}\bigr) \asymp (\log z)^{-1}\ . 
$$
Put $\psi(v)=v/\varphi(v)$. Thus 
\begin{equation}\label{eq:4.12}
G_v(\varepsilon) \asymp  \psi(v)P_v(z)^{-1}(\log z)^{-1}\ .
\end{equation} 
Inserting this into (4.10) we obtain an upper bound for $M_h(x)$, 
specifically 
\begin{equation}\label{eq:4.13}
M_h(x)\ll  \psi(h)A(x)(\log z)^{-1}\ .
\end{equation} 
where
\begin{equation}\label{eq:4.14}
A(x)= \sum_{\substack{n\le x\\(n,h)=1}}\psi(n)\theta(n)a(n)\ .
\end{equation} 
Recall that $\theta(n)$ is given by (2.1) and $a(n)$ is given by (3.1). By 
the monotonicity of the upper-bound sieve, for $n=ab$ we have 
$\theta(n) \le \theta(a)$ (simply remove the prime divisors of $b$ 
from the sifting range; cf. (6.41) of [FI1]). Hence, 
$$
A(x)\le \sum_{\substack{a\le x/y\\(a,h)=1}}\psi(a)\theta(a)\tau(a)
\sum_{b\le x/a} \psi(b)\lambda(b)\ . 
$$
The corresponding Dirichlet series 
$$
\sum_b\psi(b)\lambda(b)b^{-s} =\zeta(s)L(s,\chi)
\prod_p\Bigl(1+\frac{1}{(p-1)p^s}\bigl(1+\chi(p)-\chi(p)p^{-s}\bigr)\Bigr) 
$$
has a simple pole at $s=1$ with residue $\le 2L(1,\chi)$ so the sum 
over $b\le x/a$ is bounded by $4L(1,\chi)x/a$ 
(see more precise results in Lemma 3.2 of [FI2]) 
provided that $y\ge D^4$. Hence 
$$
A(x)\le 4L(1,\chi)x\sum_{\substack{a\le x\\(a,h)=1}}
\psi(a)\theta(a)\tau(a)a^{-1} \ll L(1,\chi)x(\log x/\log z)^{r-1}\ .
$$
To see this, first extend the range of sifting by including 
the prime factors of $h$. This extension is redundant with 
the condition $(a,h)=1$. Then, omit the condition $(a,h)=1$. Next, enlarge 
$\psi(a)\tau(a)$ to $\gamma(a)$ as in Section 3. Finally, arguing as 
for $M_h(x)$ we arrive at the above bound for $A(x)$. Inserting this 
into (4.13) we get 
\begin{equation}\label{eq:4.15}
M_h(x)\ll  \psi(h)L(1,\chi)x(\log x)^{-1}(\log x/\log z)^r\ .
\end{equation} 
Inserting this into (4.2) and then into (3.4), we conclude by (2.6) 
the following nice approximation of $S_h(x)$ by $S^*_h(x)$. 
\begin{prop}
Let $z^{72}=y =x^{1/20}$ and $D^4\le z$. Let $S_h(x)$ and $S^*_h(x)$ be given by 
(1.1) and (2.7) respectively. We have 
\begin{equation}\label{eq:4.16}
S_h(x)= S^*_h(x) + xE_h(x)
\end{equation} 
with the error term satisfying 
\begin{equation}\label{eq:4.17}
E_h(x)\ll \psi(h)L(1,\chi)\log x +(\log x)^{-2} 
\end{equation} 
for every even $h$, $0< h\le H=x(\log x)^{-24}$, 
where the implied constant is absolute. 
\end{prop}

\section {\bf Variation of $S_h(x)$ in the Shift} 

Opening the sieve functions $\theta(m)$, $\theta(n)$ in (2.7) and 
inserting (2.5) for $\Lambda^*(m)$, $\Lambda^*(n)$ we get 
$$
S_h^*(x)= \ssum_{q_1, \, q_2<y}\xi_{q_1}\xi_{q_2}\ssum_{b_1, \, b_2<y}
\nu(b_1)\nu(b_2)A_h\bigl(x;[q_1,b_1],[q_2,b_2]\bigr)\ ,
$$
with 
\begin{equation}\label{eq:5.1}
A_h(x;u,v)= \ssum_{\substack{m-n=h\\m\equiv 0(\mod u)\\n\equiv 0(\mod v)}}
\lambda'(m)\lambda'(n)\ .
\end{equation} 
Recall that $m$, $n$ run over $0< m,\,n\le x$ with $(mn,h)=1$. Note that the 
last condition implies automatically that $(q_1q_2b_1b_2, h)=1$. Now we 
need an asymptotic formula for the ``congruence sums'' $A_h(x;u,v)$ 
which holds uniformly for $0<h\le H$ and $uv\le y^4$ This a problem of 
shifted convolution type for the divisor-like function $\lambda'$, see (2.3). 
There is a vast literature on related subjects (see for example [DFI], [KMV]) 
but no result is stated which would exactly cover our sum (5.1). 
The closest seems to be Proposition 15.1 of [CI], which we adopt here without 
repeating the involved arguments (Kloosterman circle method with Weil's 
bound for Kloosterman sums). Fortunately, we do not need to use the 
results in an explicit form. In our current situation these arguments yield 
\begin{equation}\label{eq:5.2}
\begin{aligned}
& A_h(x;u,v)  =  {\mathfrak S}(h;u,v)(uv)^{-1}B(x;u,v) \\
 &\quad\quad  +O\bigl(\tau(h)(uvD)^6x^{\frac34} (\log x)^4 
+\tau(h)(uv)^{-1}x (\log x)^{-20} \bigr)\ .
\end{aligned}
\end{equation} 
The error term is negligible and $B(x;u,v)$ does not depend on $h$; 
$$
B(x;u,v) \ll x(\log x)^2 
$$
by trivial estimations. The dependence on $h$ appears 
in the ``singular series'' ${\mathfrak S}(h)={\mathfrak S}(h;u,v)$. 
The key feature of ${\mathfrak S}(h)$ is that it varies 
only slightly with respect to large prime divisors of $h$ 
(the dependence on $u$, $v$ is negligible). Specifically, 
we can write ${\mathfrak S} =1*\delta$ with
\begin{equation}\label{eq:5.3}
\delta(d)\ll \tau(d)d^{-1}\ .
\end{equation} 
Hence, for $(k,h)=1$ we derive 
\begin{equation*}
\begin{aligned}
{\mathfrak S} (hk) -{\mathfrak S} (h) 
& = \sum_{d|hk}\delta(d) -\sum_{d|h}\delta(d) =\sum_{a|h}
\sum_{\substack{c|k\\c>1}}\delta(ac) \\ 
& \ll \Bigl(\sum_{a|h}\frac{\tau(a)}{a}\Bigr) 
\Bigl(\sum_{\substack{c|k\\c>1}}\frac{\tau(c)}{c} \Bigr) 
\ll \psi(h)^2\psi(k)(\psi(k)-1)\ , 
\end{aligned}
\end{equation*} 
where we recall that $\psi(k)=k/\varphi(k)$. 
Note that for $k$ free of small prime divisors $\psi(k)-1$ is small. 

By the above observations we conclude the following 
\begin{prop} 
Assume the conditions as in Proposition 4.1. Then, for every $k$ 
with $(h,k)=1$, $1\le k\le H/h$, we have 
\begin{equation}\label{eq:5.4}
S^*_h(x) = S^*_{hk}(x) +O\bigl((\psi(k)-1)x(\log x)^6 + x(\log x)^{-2} \bigr)
\end{equation} 
where the implied constant is absolute. 
\end{prop} 

Combining (5.4) with (4.16) we obtain (under the above conditions)
\begin{equation}\label{eq:5.5}
S^*_h(x) = S^*_{hk}(x) +xE_{h,k}(x) 
\end{equation} 
where the error term satisfies 
\begin{equation}\label{eq:5.6}
E_{h,k}(x)\ll \psi(hk)L(1,\chi)\log x 
+ (\psi(k)-1)x(\log x)^6 + x(\log x)^{-2} \ .
\end{equation} 

\section {\bf Proof of Theorem 1.} 

Fix a positive even number $h$. We shall average the relation (5.5) 
with respect to $k$ over the set 
$$
{\mathcal K}= \{k:0<k\le K\ ,\, (k,hP)=1\} 
$$
where $K=H/h$ and $P=P(w)$ is the product of all primes $p\nmid h$, $p<w$ 
with $w=(\log x)^{2016}$. Note that (by elementary sieve methods) 
\begin{equation}\label{eq:6.1}
|{\mathcal K}| =K\prod_{p|hP}\bigl(1-\frac{1}{p} \bigr)  \{1+O(1/\log x)\} \ ,
\end{equation} 
\begin{equation}\label{eq:6.2}
\sum_{k\in {\mathcal K}} \psi(k) \ll |{\mathcal K}|\ ,
\end{equation}
and  
$$
\sum_{k \in {\mathcal K}} \bigl(\psi(k)-1 \bigr) 
\ll \sum_{k\in \mathcal K}\Bigl(\sum_{\substack{c|k\\ c>1}} c^{-1}\Bigr)
\le K\sum_{\substack{(c,hP)=1\\ c>1}}c^{-2} \ll K/w\ .
$$
Hence, (5.5) yields 
\begin{equation}\label{eq:6.3}
S_h(x) = \frac{1}{|{\mathcal K}|}\sum_{k\in {\mathcal K}} S_{hk}(x) 
+O\bigl(L(1,\chi)x\log x +x(\log x)^{-2} \bigr)\ .
\end{equation} 
Here we have 
$$
\sum_{k\in {\mathcal K}} S_{hk}(x)  = \sum_{\substack{n\le x\\ (n,hP)=1}}
\Lambda(n)\sum_{k \in {\mathcal K}}\Lambda (n+hk) +O(xw) 
$$
where the error term $O(xw)$ takes care of the condition $(n,hP)=1$ which is 
introduced here for technical reasons. 

Next,  we are going to execute the summatiom of $\Lambda(n+hk)$ over 
$k\in {\mathcal K}$. To this end we relax the condition $(k,hP)=1$ by 
means of upper-bound and lower-bound sieves of level $\Delta = x^{1/2016}$. 
First we get an upper bound as follows: 
\begin{equation*}
\begin{aligned}
\sum_{k\in {\mathcal K}}\Lambda(n+hk) & \le \sum_{\substack{d< \Delta\\ d|hP}}\xi_d 
\sum_{\substack{n<\ell \le n+hK\\ \ell \equiv n (\mod dh)}}\Lambda(\ell)  \\
& = hK\sum_{\substack{d< \Delta\\ d|hP}}\xi_d/\varphi(dh) +O\bigl(x(\log x)^{-A}\bigr) 
\end{aligned}
\end{equation*} 
by the Bombieri-Vinogradov theorem. Here $g(d)= \varphi(h)/\varphi(dh)$ is 
the relevant multiplicative sieve density function for which we get 
$$
\sum_{\substack{d< \Delta\\ d|hP}}\xi_d g(d) = \{1+O(1/\log x)\}
\prod_{p|hP}\bigl(1-g(p)\bigr)\ .
$$
Similarly we deal with the lower bound. Combining the two we obtain 
\begin{equation*}
\begin{aligned}
\sum_{k\in {\mathcal K}}\Lambda(n+hk) & = \{1+O(1/\log x)\}\frac{Kh}{\varphi(h)}
\prod_{p|hP}\bigl(1-g(p)\bigr)\\
& = BC(h)|{\mathcal K}|\{1+O(1/\log x)\} 
\end{aligned}
\end{equation*} 
where $B$ and $C(h)$ are given by (1.2), (1.3). Summing this over $n$ we find
\begin{equation}\label{eq:6.4}
\frac{1}{|{\mathcal K}|}\sum_{k\in {\mathcal K}}S_{hk}(x) 
=BC(h)x\{1+O(1/\log x)\}\ . 
\end{equation} 
This, together with (6.3), completes the proof of Theorem 1.

\medskip 
Department of Mathematics, University of Toronto

Toronto, Ontario M5S 2E4, Canada 

\medskip

Department of Mathematics, Rutgers University

Piscataway, NJ 08903, USA

\end{document}